\documentclass[reqno]{amsart}

\newtheorem{theorem}{Theorem}[section]
\newtheorem{lemma}[theorem]{Lemma}
\newtheorem{corollary}[theorem]{Corollary}
\newtheorem{proposition}[theorem]{Proposition}

\theoremstyle{definition}
\newtheorem{definition}[theorem]{Definition}

\newtheorem*{quest}{Question}
\theoremstyle{remark}
\newtheorem{remark}[theorem]{Remark}

\newcommand{\cB}{{\mathcal{B}}}

\newcommand{\cH}{{\mathcal{H}}}

\newcommand{\cK}{{\mathcal{K}}}

\newcommand{\fB}{{\mathfrak{B}}}

\newcommand{\fP}{{\mathfrak{P}}}
\newcommand{\bC}{{\mathbb{C}}}
\newcommand{\bN}{{\mathbb{N}}}
\newcommand{\bQ}{{\mathbb{Q}}}
\newcommand{\bR}{{\mathbb{R}}}
\newcommand{\bZ}{{\mathbb{Z}}}

\renewcommand{\Re}{{\mathrm{Re}}}

\numberwithin{equation}{section}

\newcommand{\Tr}{\mathrm{Tr}}
\newcommand{\id}{\mathrm{id}}

\newcommand{\ov}[1]{\overline{#1}}

\newcommand{\jed}{{\mathbb{I}}}
\newcommand{\spa}{\mathrm{span}}
\newcommand{\rank}{\mathrm{rank}}
\newcommand{\fA}{{\mathfrak A}}
\newcommand{\is}[2]{\langle #1 , #2\rangle}
\newcommand{\ext}{\mathrm{Ext}}

\title[Extremal positive maps]{On extremal positive maps acting between type I factors}
\author[M. Marciniak]{Marcin Marciniak}
\address{Institute of Theoretical Physics and Astrophysics, Gda{\'n}sk University,
Wita Stwosza 57, 80-952 Gda{\'n}sk, Po\-land}
\email{matmm@univ.gda.pl}
\date{11 December 2008}
\thanks{Research supported by the MNiSW research grant P03A 013 30}

\begin{document}
\setlength{\baselineskip}{1.1\baselineskip}
\begin{abstract}
The paper is devoted to the problem of classification of extremal positive maps
acting between $\fB(\cK)$ and $\fB(\cH)$ where $\cK$ and $\cH$ are Hilbert spaces.
It is shown that every positive map with the property that $\rank\,\phi(P)\leq 1$
for any one-dimensional projection $P$ is a rank 1 preserver. It allows to characterize
all decomposable extremal maps as those which satisfy the above condition.
Further, we prove that every extremal positive map which is $2$-positive turns out to 
automatically completely positive. Finally we get the same conclusion for such extremal positive maps
that $\rank\,\phi(P)\leq 1$ for some one-dimensional projection $P$ and satisfy the condition of
local complete positivity. It allows us to give a negative answer for Robertson's problem
in some special cases.
\end{abstract}
\maketitle


\section{Introduction}
Let us start with seting up some notation and terminology.
A non\-emp\-ty subset $K$ of a real or complex linear space $V$ is called
a \textit{cone} if $\alpha v +\beta w\in K$ for any $v,w\in K$ and numbers $\alpha,\beta\geq 0$.
$K$ is said to be \textit{pointed} if $K\cap(-K)=\{0\}$, and \textit{proper}
if it is pointed and closed and spans $V$. A cone $K$ in $V$ induces a partial order if we define
$v\leq w$ to mean $w-v\in K$. 
We say that a subset $F\subseteq\cK$ is a \textit{face} of $K$ if $F$ is a cone and for any $v,w\in K$ the
conditions $0\leq v\leq w$ and $w\in F$ imply $v\in F$. 
An element $v\in K$ is said to be {\it extremal} if $\{\lambda v:\,\lambda\geq 0\}$ is a face of $K$.
The set of extremal elements of $K$ we will denote by $\ext\,K$.

If $\cH$ is a Hilbert space then by $\fB(\cH)$ we denote the $C^*$-algebra of bounded operators on $\cH$.
Given a $C^*$-algebra $\fA$ and $k,l\in\bN$ we denote by 
$M_{k,l}(\fA)$ the space of all matrices of size $k\times l$ with coefficients from $\fA$. 
If $k=l$ then we will write $M_k(\fA)$ instead of $M_{k,k}(\fA)$. 
Note that $M_k(\fA)$ is canonically
isomorphic to the tensor product $M_k(\bC)\otimes\fA$, so it is endowed with the structure of a $C^*$-algebra.

Assume that $\phi:\fA\to\fB(\cH)$ is a bounded linear map. For any $k\in\bN$ we define maps 
$\phi_k,\phi^k:M_k(\fA)\to M_k(\fB(\cH))$ by
$\phi_k([A_{ij}])=[\phi(A_{ij})]$ and $\phi^k([A_{ij}])=[\phi(A_{ji})]$.
For any $C^*$-algebra $\fA$ let $\fA^+$ denote the cone of positive elements of $\fA$. 
We say that $\phi$ is a {\it positive map} whenever $\phi(\fA^+)\subseteq\fB(\cH)^+$.
We will denote the set of all positive maps from $\fA$ into $\fB(\cH)$ by $\fP(\fA,\cH)$.
If $k\in\bN$ then we say that the map $\phi$ is \textit{$k$-positive} (resp. \textit{$k$-copositive})
if the map $\phi_k$ (resp. $\phi^k$) is positive.
Whenever a map $\phi$ is $k$-positive (resp. $k$-copositive) for any $k\in\bN$ then $\phi$ is said to be 
a \textit{completely positive} (resp. \textit{completely copositive}) map.
A map $\phi$ is called \textit{decomposable} if $\phi=\phi_1+\phi_2$ for some
completely positive map $\phi_1$ and completely copositive $\phi_2$.

In spite of great efforts of many mathematicians the classification of positive maps on $C^*$-algebras 
is still a big challenge. Although there are many partial results included in several papers in mathematics as well as in mathematical physics, it
seems that we are far from full understanding of all features of these objects. For example, no algebraic
formula of a general positive map even in the case of finite dimensional matrix algebras is known.

One of the most important unsolved problems in this area is the characterization of 
extremal elements in the cone of all positive maps.
The explicit form of extremal positive unital maps is described fully only for the simplest 
non-trivial case of $2\times 2$ complex matrices (\cite{St63}). Let us warn that in our paper
we consider a larger class of all positive (i.e. not necessarily unital) maps. These both classes have
a little bit different structures. Positive unital maps form a convex subset of the cone of all unital
maps but it is not a base for this cone in the sense of \cite{AS}. As it was shown in \cite{MM06} even in
the case of $2\times 2$ matrices the structure
of extremal positive unital maps differs from the structure of extremal elements in the cone of all positive maps.

On the other hand, let us remind that all extremal elements of the cone of completely positive maps
are fully recognized (see \cite{Choi75a,Arv69}).
If we consider maps from $\fB(\cK)$ into $\fB(\cH)$ where $\cK$ and $\cH$ are finite dimensional Hilbert spaces, 
then a map $\phi$ is extremal in the cone of completely
positive maps if and only if $\phi(X)=AXA^*$, $X\in\fB(\cK)$, where $A\in\fB(\cK,\cH)$. Analogously, extremal
maps in the cone of completely copositive maps are of the form $\phi(X)=AX^{\mathrm T}A^*$, $X\in\fB(\cK)$, for some
$A\in\fB(\cK,\cH)$, where $X^{\mathrm T}$ denotes a transposition of the element $X$.
Consequently, the cone of all decomposable maps
is the hull of all maps which have one of the previously
mentioned two forms.
Coming back to positive maps, it was proved in \cite{YH05} that  maps of the above two forms are extremal also in the cone $\fP(\fB(\cK),\cH)$ of all positive maps.

Further, let us note that in cones $\fP(\fB(\bC^2),\bC^2)$, $\fP(\fB(\bC^2),\bC^3)$ and $\fP(\fB(\bC^3),\bC^2)$ there are no other than the mentioned above extremal elements. This is a consequence of the results of St{\o}rmer (\cite{St63})
and Woronowicz (\cite{Wor76}) that these cones contain only decomposable maps.
However, there are known some other examples of extremal positive maps between matrix algebras for greater dimensions
(see \cite{Choi75b,Ha03,Kye95}). 
Obviously, they are necessarily nondecomposable. The most famous example is that which belongs to $\fP(\fB(\bC^3),\bC^3)$ given by Choi in \cite{Choi75b}
\begin{equation}\label{Choi}
\phi\left(\left[\begin{array}{ccc}
a_{11}&a_{12}&a_{13}\\
a_{21}&a_{22}&a_{23}\\
a_{31}&a_{32}&a_{33}
\end{array}\right]\right)=
\left[\begin{array}{ccc}
a_{11}+a_{33}&-a_{12}&-a_{13}\\
-a_{21}&a_{22}+a_{11}&-a_{23}\\
-a_{31}&-a_{32}&a_{33}+a_{22}
\end{array}\right] .
\end{equation}
It was the first known example of a nondecomposable map. Let us mention also that in the literature 
several examples of nondecomposable maps are described (\cite{ChKyLe,ChKo06,Ha03,Kye95,Osa91,Rob83,St80a,Tang86,Ter01}).
Although some conditions equivalent to decomposability are known (\cite{St80}), proving that a positive
map is nondecomposable is a very difficult task. But it seems that providing new examples of extremal
maps is of extremal difficulty.

Apart from the mentioned above results there is another line of research in the mathematical literature 
which deals with similar problems. It comes from convex analysis (see \cite{Tam95} and references therein).
The main object in this framework is an ordered linear space, i.e. a pair $(V,V^+)$ where $V$ is a finite dimensional 
linear
space while $V^+$ is a pointed cone in $V$. Having two such objects, say $(V,V^+)$ and $(W,W^+)$ we can
consider maps $T:V\to W$ such that $T(V^+)\subset W^+$. We call them {\it positive maps}, and they
form a cone which we will denote by $\fP(V,W)$.
As in the case of $C^*$-algebras we try to describe extremal elements of that cone. 
Let us remind an interesting result of Loewy and Schneider which goes in this direction.
To this end we recall that a cone $K$ is {\it indecomposable} if there are no non-empty subsets $K_1,K_2\subset K$
such that $K=K_1+K_2$ and $\spa\, K_1\cap\spa\, K_2=\{0\}$. 
\begin{theorem}[\cite{LS75}]
Let $V^+$ be a cone in $V$ and assume $V^+=\mathrm{hull}(\ext\,V^+)$. Then the following conditions
are equivalent:
\begin{enumerate}
\item[(1)] $V^+$ is indecomposable.
\item[(2)] If $T\in L(V,V)$ is such that $\ker T=\{0\}$ and $T(\ext\,V^+)\subseteq\ext\,V^+$,
then $T\in\ext\fP(V,V)$.
\item[(3)] If $T\in L(V,V)$ is such that $\ker T=\{0\}$ and $T(V^+)=V^+$, then $T\in\ext\fP(V,V)$.
\item[(4)] $\id_V\in\ext\fP(V,V)$.
\end{enumerate}
\end{theorem}
Now, let $\cH$ be a finite dimensional Hilbert space, $V=\fB(\cH)_{\mathrm h}$ be the space
of selfadjoint elements of $\fB(\cH)$, and $V^+=\fB(\cH)^+$. It was proved in \cite{YH05} that $\fB(\cH)^+$ is an indecomposable cone. Note also that $\ext\fB(\cH)^+$ consists of nonnegative multiplicities of one-dimensional projections on $\cH$. Hence, we conclude that if $\phi:\fB(\cH)\to\fB(\cH)$
is a bijective linear mapping such that $\rank\,\phi(P)=1$ for any one-dimensional projection $P$ then
$\phi$ is extremal in $\fP(\fB(\cH),\cH)$. 
\begin{remark}\label{rem1}
We will see later that the implication contained in the point (2) of
the above theorem can not be converted. The Choi map (\ref{Choi}) will serve as a counterexample, because it sends
all one-dimensional projections into operators of rank not smaller than $2$. But it is still an open problem
whether there exists $T\in\ext\fP(V,V)$ such that $T(\ext\,V^+)\not\subseteq\ext\,V^+$ but $T(\ext\,V^+)\cap\ext\,V^+\neq\emptyset$. In the context of operator algebras we will prove in Section 2
that such a map must be nondecomposable. So, we can ask if there is
a nondecomposable extremal positive map such that $\rank\,\phi(P)=1$ for some one-dimensional projection $P$.
\end{remark}

The last statement before the above remark tell us to draw our attention to the theory 
of the so called linear preservers (for survey see \cite{LinPres}).
In particular, we are interested in the problem of rank $1$ preservers i.e. 
linear maps $T:\fB(\cK)\to\fB(\cH)$
such that $\rank\,\phi(X)=1$ whenever $\rank\,X=1$ for $X\in\fB(\cK)$.
They are well described. By the result of Marcus and Moyls (\cite{MarMoy59})
we know that each injective rank $1$ preserver is of the form $T(X)=MXN$ or 
$T(X)=MX^{\mathrm T}N$, $X\in\fP(\cK)$,
for some $M\in\fB(\cK,\cH)$ and $N\in\fB(\cH,\cK)$. 
Lim (\cite{Lim75}) proved that the similar form follows from a weaker assumption. Namely,
it is enough to assume that $\rank\,\phi(X)\leq 1$ for every $X\in\fB(\cK)$ such that $\rank\,X=1$.
\begin{remark}\label{rem2}
Observe that if a map $\phi\in\fP(\fB(\cK),\cH)$ is of the form $\phi(X)=AXA^*$ or $\phi(X)=AX^{\mathrm T}A^*$
then $\rank\,\phi(P)\leq 1$ for every one-dimensional projection $P$ on $\cK$.
Motivated by the considerations from the above paragraph we can ask whether the converse is true (cf. \cite{GMS98} 
and references therein).
\end{remark}

The aim of this paper is to present a new approach to the problem of classification
of extremal maps in $\fP(\fB(\cK),\cH)$ which is based on point of view coming from 
the convex analysis and linear preservers theory. Our main motivation is to give answers 
for questions contained in Remarks \ref{rem1} and \ref{rem2}.
However, if it comes to our methods we will use the technique presented in papers 
\cite{LMM06,MM01,MM05,MM06,MM07}.

The paper is organized as follows. In Section 2 we give an 'almost' positive answer for the question from
Remark \ref{rem2} (Theorem \ref{t:rank1}).
It will allow us to characterize decomposable extremal maps as those maps which have rank $1$ 
nonincreasing property (Corollary \ref{m}).
In Section 3 we formulate some conditions on a map $\phi$ which are equivalent to the property that
$\phi$ is minorized by some completely positive (or completely copositive) extremal map (Theorem \ref{p:psil}).
As a consequence we get the result that each extremal map in $\fP(\fB(\cK),\cH)$ which is $2$-positive (resp.
$2$-copositive) is automatically completely positive (resp. completely copositive) (Theorem \ref{t:2pos}). 
It is a partial negative answer
for the question asked by Robertson in \cite{Rob83}.
The aim of Section 4 is to show that under some continuity assumptions each positive map can be reconstructed
from its values on one-dimensional projections (Theorem \ref{t:par}).
In the last section we deal with the problem formulated in Remark \ref{rem1}.
Firstly, motivated by the results of Section 4, we describe properties of such maps $\phi$  
that $\rank\,\phi(P)=1$ for some one-dimensional projection $P$ in terms of some positive functions with
parallelogram identity.
We apply the technique developed in our previous papers (\cite{MM06,MM07}) to consider 
such maps which are extremal. We show in Theorem \ref{t:G-F} that under some additional relatively
weak condition of local complete positivity they are completely positive. This is a partial negative solution
for Robertson's problem as well as a strong suggestion that the problem from Remark \ref{rem1}
has negative solution. Finally, in Corollary \ref{t:G-F2} we give negative answer for Robertson's question 
in the case when $\cH$ is finite dimensional and $\dim\cK=2$.

\section{Rank $1$ nonincreasing positive maps}
Let $\cK$ and $\cK$ be Hilbert spaces. For now we do not formulate any assumptions about the dimensions of these
spaces, but in the sequel it may happen that we will need to work with finite dimensional spaces.

Let us introduce some notations. 
If $\xi,\eta\in\cK$ then by $\xi\eta^*$ we denote the operator on $\cK$ which is defined by
$$(\xi\eta^*)\tau=\langle\eta,\tau\rangle\xi,\qquad\tau\in\cK.$$
We assume also that some antilinear selfadjoint involution on $\cK$ is defined, i.e. such a map
$\cK\ni\xi\mapsto\overline{\xi}\in\cK$ that
\begin{enumerate}
\item $\overline{a\xi+b\eta}=\overline{a}\overline{\xi}+\overline{b}\overline{\eta}$ for $\xi,\eta\in\cK$ and $a,b\in\bC$,
\item $\overline{\overline{\xi}}=\xi$ for $\xi\in\cK$,
\item $\langle\overline{\xi},\overline{\eta}\rangle=\langle\eta,\xi\rangle$ for $\xi,\eta\in\cK$.
\end{enumerate}
Having such an involution one can define for every $X\in\fB(\cK)$ its {\em transpose} $X^{\mathrm T}$ (with respect to the
involution) by
$$X^{\mathrm T}\xi=\overline{{X^*} \overline{\xi}},\qquad\xi\in\cK.$$
We observe that the transposition is a linear $^*$-antimorphism. Moreover,
$(\xi\eta^*)^{\mathrm T}=\overline{\eta}\overline{\xi}^*$ for any $\xi,\eta\in\cK$.

The main task of this section is to describe all positive maps $\phi:\fB(\cK)\to \fB(\cH)$ which have the 
property that $\rank\, \phi(P)\leq 1$ for every $1$-dimensional projection $P$ acting on $\cK$. Such maps we will call
\textit{rank-$1$ nonincreasing positive maps}. We start with the following
\begin{lemma}\label{lemma}
Let $x,y\in\cH$. Assume that $A\in \fB(\cH)$ satisfies 
$$\rank (xx^*+|\lambda|^2yy^*+\lambda A+\ov{\lambda}A^*)\leq 1$$
for every $\lambda\in\bC$. Then 
\begin{enumerate}
\item\label{lemindep} $A=\mu xy^*$ or $A=\mu yx^*$ for some $\mu\in\bC$ with
$|\mu|=1$ whenever $x$ and $y$ are linearly independent;
\item\label{lemdepx} $A=\mu xx^*$ for some $\mu\in\bC$ when $x\neq 0$ and $x,y$ are linearly dependent;
\item\label{lemdepy} $A=\mu yy^*$ for some $\mu\in\bC$ when $y\neq 0$ and $x,y$ are linearly dependent;
\item\label{lemzero} $A$ is a complex multiplicity of some one-dimensional projection in $\fB(\cH)$ if $x=y=0$.
\end{enumerate}
\end{lemma}
\begin{proof}
Let us denote 
$R_\lambda=xx^*+|\lambda|^2yy^*+\lambda A + \ov{\lambda}A^*$ for any $\lambda\in\bC$.
Let $\xi$ and $\eta$ be any vectors from $\cH$.
By the assumption $R_\lambda$ has rank at most one, so the vectors $R_\lambda \xi$ and $R_\lambda \eta$ are linearly dependent. 
Hence $\delta=0$ where
$\delta=\is{\xi}{R_\lambda\xi}\is{\eta}{R_\lambda\eta}-|\is{\xi}{R_\lambda\eta}|^2$.
We can calculate that
$\delta 
=\alpha_1(\theta)r+\alpha_2(\theta)r^2+\alpha_3(\theta)r^3$
where $r\geq 0$ and $\theta\in[0,2\pi)$ are such that $\lambda=re^{i\theta}$, and
\begin{eqnarray}
\alpha_1(\theta)&=&
2|\is{\eta}{x}|^2\Re\, e^{i\theta}\is{\xi}{A\xi} + 2|\is{\xi}{x}|^2\Re\, e^{i\theta}\is{\eta}{A\eta}  \label{alfa1}\\
&&-
2\Re\, e^{i\theta}\left(\is{x}{\xi}\is{\xi}{A\eta}\is{\eta}{x}+\is{x}{\eta}\is{\eta}{A\xi}\is{\xi}{x}\right), \nonumber\\[2mm]
\alpha_2(\theta)&=& \left|\is{\xi}{x}\is{\eta}{y}-\is{\xi}{y}\is{\eta}{x}\right|^2 + 
4\left(\Re\, e^{i\theta}\is{\xi}{A\xi}\right)\left(\Re\, e^{i\theta}\is{\eta}{A\eta}\right) \label{alfa2}\\
&&-\left|e^{i\theta}\is{\xi}{A\eta}+e^{-i\theta}\ov{\is{\eta}{A\xi}}\right|^2, \nonumber\\[2mm]
\alpha_3(\theta)&=&
2|\is{\eta}{y}|^2\Re\, e^{i\theta}\is{\xi}{A\xi} + 2|\is{\xi}{y}|^2\Re\, e^{i\theta}\is{\eta}{A\eta} \label{alfa3}\\
&&- 
2\Re\, e^{i\theta}\left(\is{y}{\xi}\is{\xi}{A\eta}\is{\eta}{y}+\is{y}{\eta}\is{\eta}{A\xi}\is{\xi}{y}\right) .\nonumber
\end{eqnarray}
Let us fix $\theta$ for a moment. Then $\delta$ becomes a polynomial of the real variable $r$.
Since it is zero for any $r>0$, each of the three coefficients of this polynomial should vanish for any $\theta\in[0,2\pi)$.

In order to prove the statement (a) we assume that the vectors $x$ and $y$ are linearly independent.
Let us consider the following three special cases of the choice of $\xi$ and $\eta$:

\textit{1st case:} $\xi=x$ and $\eta$ orthogonal to both $x$ and $y$.
Then formulas (\ref{alfa1}) and (\ref{alfa2}) reduce to
\begin{eqnarray*}
\alpha_1(\theta)&=&
2\Vert x\Vert^4\Re\, e^{i\theta}\is{\eta}{A\eta},\\
\alpha_2(\theta)&=&
4\left(\Re\, e^{i\theta}\is{x}{Ax}\right)\left(\Re\, e^{i\theta}\is{\eta}{A\eta}\right)
-\left|e^{i\theta}\is{x}{A\eta}+e^{-i\theta}\ov{\is{\eta}{Ax}}\right|^2, 
\end{eqnarray*}
Since $\alpha_1(\theta)=0$ for any $\theta$, we infer that 
\begin{equation}\label{etaAeta}
\is{\eta}{A\eta}=0.
\end{equation}
Next, $\alpha_2(\theta)=0$ implies that 
$e^{i\theta}\is{x}{A\eta}+e^{-i\theta}\ov{\is{\eta}{Ax}}=0$
for any $\theta$, and consequently
\begin{equation}\label{xAeta}
\is{x}{A\eta}=0\qquad\textrm{and}\qquad\is{\eta}{Ax}=0.
\end{equation}

\textit{2nd case:} $\xi=y$ and $\eta$ orthogonal to both $x$ and $y$.
By similar arguments as in the previous case we obtain
\begin{equation}\label{yAeta}
\is{y}{A\eta}=0\qquad\textrm{and}\qquad\is{\eta}{Ay}=0.
\end{equation}

\textit{3rd case:} $\xi=x$ and $\eta=y$. In this case the formulas 
(\ref{alfa1}), (\ref{alfa2}) and (\ref{alfa3})
take the form
\begin{eqnarray}
\alpha_1(\theta)&=&
2\left(|\is{x}{y}|^2\Re\, e^{i\theta}\is{x}{Ax} + \Vert x\Vert^4\Re\, e^{i\theta}\is{y}{Ay}\right) \nonumber\\
&&-2 
\Vert x\Vert^2\Re\, e^{i\theta}(\is{y}{x}\is{x}{Ay}+\is{x}{y}\is{y}{Ax}), \nonumber\\
\alpha_2(\theta)&=& \left(\Vert x\Vert^2\Vert y\Vert^2-|\is{x}{y}|^2\right)^2 + 
4\left(\Re\, e^{i\theta}\is{x}{Ax}\right)\left(\Re\, e^{i\theta}\is{y}{Ay}\right) \label{e:a2}\\
&&-\left|e^{i\theta}\is{x}{Ay}+e^{-i\theta}\ov{\is{y}{Ax}}\right|^2, \nonumber\\
\alpha_3(\theta)&=&2\left(\Vert y\Vert^4\Re\, e^{i\theta}\is{x}{Ax} + |\is{x}{y}|^2\Re\, e^{i\theta}\is{y}{Ay}\right) \nonumber\\
&&- 2
\Vert y\Vert^2\Re e^{i\theta}(\is{y}{x}\is{x}{Ay}+\is{x}{y}\is{y}{Ax}) .\nonumber
\end{eqnarray}

The equalities $\alpha_1(\theta)=0$ and $\alpha_3(\theta)=0$ imply the following conditions
\begin{eqnarray}
\Re\, e^{i\theta}\is{x}{Ax}&=&\frac{\Vert x\Vert^2\Re\, e^{i\theta}
\left(\is{y}{x}\is{x}{Ay}+\is{x}{y}\is{y}{Ax}\right)}{\Vert x\Vert^2\Vert y\Vert^2+|\is{x}{y}|^2}, \label{RexAx}\\
\Re\, e^{i\theta}\is{y}{Ay}&=&\frac{\Vert y\Vert^2\Re\, e^{i\theta}
\left(\is{y}{x}\is{x}{Ay}+\is{x}{y}\is{y}{Ax}\right)}{\Vert x\Vert^2\Vert y\Vert^2+|\is{x}{y}|^2}. \label{ReyAy}
\end{eqnarray}
If we substitute both these expressions into (\ref{e:a2}) 
then it turns out that the
equation $\alpha_2(\theta)=0$ is equivalent to the following
\begin{eqnarray*}
\lefteqn{\left(\Vert x\Vert^2\Vert y\Vert^2-|\is{x}{y}|^2\right)^2-
\left(|\is{x}{Ay}|^2+|\is{y}{Ax}|^2\right)=}\\
&=&2\Re\, e^{2i\theta}\is{x}{Ay}\is{y}{Ax} -
\frac{4\Vert x\Vert^2\Vert y\Vert^2\left[\Re\, e^{i\theta}(\is{y}{x}\is{x}{Ay}+\is{x}{y}\is{y}{Ax})\right]^2}{\left(\Vert x\Vert^2\Vert y\Vert^2+|\is{x}{y}|^2\right)^2}
.
\end{eqnarray*} From 
the identity $(\Re z)^2=\frac{1}{2}\Re z^2+\frac{1}{2}|z|^2$ for any $z\in \bC$ we infer that the above 
equality is equivalent to 
\begin{equation}\label{bety}
\beta_1=2\Re\, e^{2i\theta}\beta_2,
\end{equation}
where
\begin{eqnarray}
\beta_1&=&\left(\Vert x\Vert^4\Vert y\Vert^4-|\is{x}{y}|^4\right)^2 \\ &&-\;\left(|\is{x}{Ay}|^2+|\is{y}{Ax}|^2\right)
\left(\Vert x\Vert^4\Vert y\Vert^4+|\is{x}{y}|^4\right) \nonumber\\
&&+\;4\Vert x\Vert^2\Vert y\Vert^2\Re\is{y}{x}^2\is{x}{Ay}\ov{\is{y}{Ax}} \nonumber\\[2mm]
\beta_2&=&
\left(\Vert x\Vert^2\Vert y\Vert^2\is{y}{Ax}-\is{y}{x}^2\is{x}{Ay}\right)\times\\ &&\times\;
\left(\Vert x\Vert^2\Vert y\Vert^2\is{x}{Ay}-\is{x}{y}^2\is{y}{Ax}\right)\nonumber
\end{eqnarray}
Let us observe that the condition (\ref{bety}) holds for any $\theta$. It is possible if and only if
$\beta_1=0$ and $\beta_2=0$. The last equality implies that
\begin{equation}\label{beta2a}
\Vert x\Vert^2\Vert y\Vert^2\is{y}{Ax}=\is{y}{x}^2\is{x}{Ay}
\end{equation}
or
\begin{equation}\label{beta2b}
\Vert x\Vert^2\Vert y\Vert^2\is{x}{Ay}=\is{x}{y}^2\is{y}{Ax}
\end{equation}
Assume that (\ref{beta2a}) is satisfied. Then
$$\is{y}{Ax}=\frac{\is{y}{x}^2}{\Vert x\Vert^2\Vert y\Vert^2}\is{x}{Ay}.$$
Considering the fact that $\beta_1=0$ this leads to the equality
$$\left(\Vert x\Vert^4\Vert y\Vert^4-|\is{x}{y}|^4\right)^2\left(\Vert x\Vert^4\Vert y\Vert^4-|\is{x}{Ay}|^2\right)=0.$$
which implies
\begin{equation}\label{mod1}
|\is{x}{Ay}|=\Vert x\Vert^2\Vert y\Vert^2.
\end{equation}
As a consequence of (\ref{beta2a}) and (\ref{mod1}) we get
\begin{eqnarray}\label{xAy}
\is{x}{Ay}&=&\mu\Vert x\Vert^2\Vert y\Vert^2 \\
\is{y}{Ax}&=&\mu\is{y}{x}^2 \label{yAx}
\end{eqnarray}
for some complex number $\mu$ such that $|\mu|=1$. 

Since (\ref{RexAx}) and (\ref{ReyAy}) hold for any $\theta$, the following conditions must be satisfied
\begin{eqnarray}
\is{x}{Ax}&=&\frac{\Vert x\Vert^2}{\Vert x\Vert^2\Vert y\Vert^2+|\is{x}{y}|^2}
(\is{y}{x}\is{x}{Ay}+\is{x}{y}\is{y}{Ax}) \label{pre_xAx}\\
\is{y}{Ay}&=&\frac{\Vert y\Vert^2}{\Vert x\Vert^2\Vert y\Vert^2+|\is{x}{y}|^2}
(\is{y}{x}\is{x}{Ay}+\is{x}{y}\is{y}{Ax}) \label{pre_yAy}
\end{eqnarray}
If we apply (\ref{xAy}) and (\ref{yAx}) into (\ref{pre_xAx}) and (\ref{pre_yAy}) then we obtain
\begin{eqnarray}
\is{x}{Ax}&=&\mu\Vert x\Vert^2\is{y}{x}, \label{xAx}\\
\is{y}{Ay}&=&\mu\Vert y\Vert^2\is{y}{x}. \label{yAy}
\end{eqnarray}

In the same way one can deduce from (\ref{beta2b})
the following set of relations
\begin{eqnarray}
\label{xA'y}
\is{x}{Ay}&=&\mu\is{x}{y}^2 \\
\is{y}{Ax}&=&\mu\Vert x\Vert^2\Vert y\Vert^2, \label{yA'x} \\
\is{x}{Ax}&=&\mu\Vert x\Vert^2\is{x}{y}, \label{xA'x}\\
\is{y}{Ay}&=&\mu\Vert y\Vert^2\is{x}{y}. \label{yA'y}
\end{eqnarray}

Let us summarize the results contained in the above three cases. 
Denote by $P$ the orthogonal projection onto the subspace of $\cH$ 
generated by $x$ and $y$. 
Since (\ref{etaAeta}) holds for any $\eta$ orthogonal to the subspace $P\cH$, it follows that
\begin{equation}\label{QAQ}
(\jed-P)A(\jed-P)=0.
\end{equation}
Further, from the fact that (\ref{xAeta}) and (\ref{yAeta}) are satisfied for any $\eta$ orthogonal to
$P\cH$ we conclude that
\begin{equation}\label{QAP}
(\jed-P)AP=0\qquad\textrm{and}\qquad PA(\jed-P)=0.
\end{equation}
Finally, we discovered in the 3rd case that $A$ fulfils one the following two sets of relations:
(\ref{xAy}), (\ref{yAx}), (\ref{xAx}), (\ref{yAy}) or 
(\ref{xA'y}), (\ref{yA'x}), (\ref{xA'x}), (\ref{yA'y}). It follows from independence of $x$ and $y$
that
\begin{equation}\label{PAP}
PAP=\mu xy^*\qquad\textrm{or}\qquad PAP=\mu yx^*
\end{equation}
for some complex number $\mu$ such that $|\mu|=1$.
We finish the proof of the statement (1) by the observation that it follows from (\ref{QAQ}), (\ref{QAP}) and (\ref{PAP}).

Now, let us assume that $x$ and $y$ are linearly dependent.   
Then, it is easy to observe that for any $\xi,\eta\in\cH$ and $\theta\in[0,2\pi)$ 
the formula
(\ref{alfa2}) has the form
\begin{equation}
\alpha_2(\theta)=  
4\left(\Re\, e^{i\theta}\is{\xi}{A\xi}\right)\left(\Re\, e^{i\theta}\is{\eta}{A\eta}\right)
-\left|e^{i\theta}\is{\xi}{A\eta}+e^{-i\theta}\ov{\is{\eta}{A\xi}}\right|^2, \label{alfa2d} 
\end{equation}
If $A=0$ then it satisfies each of the statements (b), (c) and (d), so without loss of generality we may assume $A\neq 0$. Then there exists $\xi_0\in\cH$ such that $\is{\xi_0}{A\xi_0}\neq 0$. Let 
\begin{equation}\label{polar}
\is{\xi_0}{A\xi_0}=a e^{it}
\end{equation}
for some $a>0$ and $t\in[0,2\pi)$. Let $\eta$ be an arbitrary vector and $\xi=\xi_0$.
As $\alpha_2(\theta)=0$ for any $\theta$, we observe
$$\left(\Re\, e^{i\theta}\is{\xi_0}{A\xi_0}\right)\left(\Re\, e^{i\theta}\is{\eta}{A\eta}\right)
=\frac{1}{4}\left|e^{i\theta}\is{\xi_0}{A\eta}+e^{-i\theta}\ov{\is{\eta}{A\xi_0}}\right|^2\geq 0$$
So, for any $\theta$ both real numbers 
$\Re\, e^{i\theta}\is{\xi_0}{A\xi_0}$ and $\Re\, e^{i\theta}\is{\eta}{A\eta}$ are of the same sign. It is possible 
in the case when both complex numbers 
$\is{\xi_0}{A\xi_0}$ and $\is{\eta}{A\eta}$ have the same argument, i.e. 
\begin{equation}\label{polar1}
\is{\eta}{A\eta}=b_\eta e^{it}
\end{equation}
for some $b_\eta\geq 0$ and the argument $t$ is determined by (\ref{polar}).
Let $B=e^{-it}A$. Then it follows from (\ref{polar1}) that for any $\eta\in\cH$ we have
$\is{\eta}{B\eta}=b_\eta\geq 0$, hence $B$ is a positive operator.

In order to prove (b) let us consider $x\neq 0$ and $y=\gamma x$ for some $\gamma\in \bC$.
In the same way (c.f. the paragraph containing formula (\ref{etaAeta})) as in the 1st case in the proof of statement (1) we prove that
$\is{\eta}{A\eta}=0$ for any $\eta$ orthogonal to $x$. This implies that $B=cxx^*$
for some $c>0$, and consequently $A=ce^{it}xx^*$.

The statement (c) follows from the similar line of arguments as above.

Now, assume $x=y=0$, so $\lambda A+\ov{\lambda} A^*$ has rank at most one for any $\lambda$.
But $\lambda A+\ov{\lambda} A^*=\left(2\Re \lambda e^{it}\right)B$, so $B$ must be a multiplicity
of a one-dimensional projection. This finishes the proof of (d).
\end{proof}

Now, we are ready to formulate the main result of this section
\begin{theorem}\label{t:rank1}
Assume that $\cK$ and $\cH$ are finite dimensional Hilbert spaces and $\phi:\fB(\cK)\to\fB(\cH)$ is a rank $1$
non-increasing positive map. Then one of the following three conditions holds:
\begin{enumerate}
\item[(i)] there exist a positive functional $\omega$ on $\fB(\cK)$ and a one-dimensional projection $Q$ on $\cH$ such that $\phi(X)=\omega(X)Q$ for any $X\in\fB(\cK)$;
\item[(ii)] there exists a linear operator $B\in\fB(\cK,\cH)$ such that $\phi(X)=BXB^*$ for any $X\in\fB(\cK)$;
\item[(iii)] there exists a linear operator $C\in\fB(\cK,\cH)$ such that $\phi(X)=CX^{\mathrm T}C^*$ for any $X\in\fB(\cK)$.
\end{enumerate}
\end{theorem}
\begin{proof}
By the assumption $\phi$ maps one-dimensional orthogonal projections into positive multiplicities
of one-dimensional projections. Hence, for any $\xi\in\cK$ there is a vector $x_\xi$ (not uniquely determined)
such that $\phi(\xi\xi^*)=x_\xi x_\xi^*$.
We will prove that $\rank\,\phi(X)\leq 1$ whenever $\rank\,X\leq 1$ for any $X\in\fB(\cK)$. 
Let $X\in\fB(\cK)$ be of rank $1$.
Then $X=\xi\eta^*$ for some vectors $\xi,\eta\in\cK$. Let $A=\phi(\xi\eta^*)$. We must show 
that $\rank\,A\leq 1$. To this end define for $\lambda\in\bC$ 
$$R_\lambda=x_\xi x_\xi^*+|\lambda|^2 x_\eta x_\eta^*+\lambda A+\ov{\lambda}A^*.$$
Observe that
$R_\lambda=\phi\left((\xi+\lambda \eta)(\xi+\lambda \eta)^*\right)$, so by assumption $\rank\, R_\lambda\leq 1$
for every $\lambda\in\bC$.
Now, we conclude from Lemma \ref{lemma} that $\rank\, A\leq 1$. 

From Theorem 1 in \cite{Lim75} (see also \cite{MarMoy59}) we conclude that one of the following conditions must
hold:
\begin{enumerate}
\item[(a)] $\phi(\fB(\cK))\setminus\{0\}$ consists entirely of rank one operators;
\item[(b)] there exist linear operators $B_1\in\fB(\cK,\cH)$, $B_2\in\fB(\cH,\cK)$ such that $\phi(X)=B_1XB_2$
for all $X\in\fB(\cK)$;
\item[(c)] there exist linear operators $C_1\in\fB(\cK,\cH)$, $C_2\in\fB(\cH,\cK)$ such that $\phi(X)=C_1X^{\mathrm T}C_2$
for all $X\in\fB(\cK)$.
\end{enumerate}

Assume that (a) is valid. Since $\phi$ is positive $\phi(\jed)=\mu Q$ for some $\mu>0$ and 
a one-dimensional projection $Q$. For any element $X\in\fB(\cK)$ the inequality $X\leq \Vert X\Vert \jed$ holds, so
$\phi(X)\leq \mu\Vert X\Vert Q$. Thus, for any $X$ there is a number $\omega(X)$ such that $\phi(X)=\omega(X)Q$.
Linearity and positivity of $\phi$ implies the same properties for $\omega$. So Condition (i) holds.

Now, assume that (b) is fulfiled. Then 
$0\leq\phi(\xi\xi^*)=B_1\xi\xi^*B_2=(B_1\xi)(B_2^*\xi)^*$ for any $\xi\in\cK$. It follows that $B_2^*\xi=\lambda_\xi B_1\xi$ for
some $\lambda_\xi\geq 0$. But for any $\xi,\eta\in\cK$ we have
\begin{eqnarray*}
\lefteqn{\lambda_\eta(B_1\xi)(B_1\eta)^*=(B_1\xi)(B_2^*\eta)^*=
\phi(\xi\eta^*)=}\\ &=&\phi(\eta\xi^*)^*= ((B_1\eta)(B_2^*\xi)^*)^*=
\lambda_\xi((B_1\eta)(B_1\xi)^*)^*=\lambda_\xi(B_1\xi)(B_1\xi)^*.
\end{eqnarray*}
Thus $\lambda_\eta=\lambda_\xi$, so $B_2=\lambda B_1^*$ for some constant $\lambda\geq 0$. If $B=\lambda^{1/2}B_1$
then $\phi(X)=BXB^*$ for any $X\in\fB(\cK)$, so (ii) is valid.

By similar arguments we show that the property (c) implies (iii).
\end{proof}

The result of \cite{YH05} asserts that maps of the form $\phi(X)=BXB^*$ and $\phi(X)=CX^{\mathrm T}C^*$ are extremal in the cone 
of positive maps between $\fB(\cK)$ and $\fB(\cH)$ provided that $\cK$ and $\cH$ are finite dimensional. By the above theorem
it is possible to characterize these maps among all extremals in terms of rank properties.
\begin{corollary}\label{m}
Assume that $\phi:\fB(\cK)\to \fB(\cH)$ is an extremal positive linear map. Then the following
conditions are equivalent:
\begin{enumerate}
\item[(1)] $\phi$ is decomposable;
\item[(2)] $\phi$ is either completely positive or completely copositive;
\item[(3)] there is $B\in\fB(\cK,\cH)$ such that $\phi(X)=BXB^*$ for all $X\in\fB(\cK)$ or
there is $C\in\fB(\cK,\cH)$ such that $\phi(X)=CX^{\mathrm T}C^*$ for any $X\in \fB(\cK)$;
\item[(4)] for any one-dimensional projection $P$ on $\cK$ there is a one-dimensional projection $Q$ and a non-negative constant $\lambda$ such that $\phi(P)=\lambda Q$. 
\end{enumerate}
\end{corollary} 
\begin{proof}
(1) $\Rightarrow$ (2) 
By the assumption $\phi=\phi_1+\phi_2$ where $\phi_1$ is a completely positive map while $\phi_2$ is
a completely copositive one. Then we have $\phi_1\leq\phi$ and $\phi_2\leq\phi$. The extremality
of $\phi$ implies $\phi_1=\lambda_1\phi$ and $\phi_2=\lambda_2\phi$ for some $\lambda_1,\lambda_2\geq 0$. 
Since $\phi\neq 0$, we have $\lambda_1>0$ or $\lambda_2>0$, and the assertion (2) is proved.

(2) $\Rightarrow$ (3)
Assume that $\phi$ is a completely positive map. Then, by the result of Choi (cf. \cite[Theorem 1]{Choi75a})
$\phi$ is of the form $\phi(X)=\sum_i^kA_iXA_i^*$ for all $X\in\fB(\cK)$ where $k\in\bN$ and $A_1,\ldots,
A_k\in\fB(\cK,\cH)$.
It follows from the extremality of $\phi$ that the sum must reduce to a one term, so $\phi$ is of the form
$\phi(X)=AXA^*$.
Suppose now that $\phi$ is completely copositive. It is equivalent to the fact that the map $X\mapsto \phi(X^{\mathrm T})$
is completely positive. So, it follows from the theorem of Choi that $\phi$ is of the form
$\phi(X)=\sum_i^kA_iX^{\mathrm T}A_i^*$, and we use the same argument based on extremality of $\phi$ as above to deduce
that $\phi(X)=AX^{\mathrm T}A^*$ for some $A\in \fB(\cK,\cH)$.

(3) $\Rightarrow$ (1) 
Obvious.

(3) $\Rightarrow$ (4)
Let $P=\xi\xi^*$ for some $\xi\in\cK$ such that $\Vert\xi\Vert=1$.
If $\phi(X)=BXB^*$ for $X\in\fB(\cK$ then $\phi(\xi\xi^*)=B\xi\xi^*B^*=(Bx)(Bx)^*$,
so $\phi(P)$ is a multiplicity of some one-dimensional projection. In the case when $\phi(X)=CX^TC^*$ for $X\in\fB(\cK)$ we have 
$\phi(\xi\xi^*)=C(\xi\xi^*)^TC^*=C\ov{\xi}\ov{\xi}^*C^*=(C\ov{\xi})(C\ov{\xi})^*$, so we get the same conclusion.

(4) $\Rightarrow$ (3)
By Theorem \ref{t:rank1} $\phi(X)=\omega(X)Q$ or $\phi(X)=BXB^*$ or $\phi(X)=CX^TC^*$. In the last two possibilities we have (3). 
Assume the first possibility. Extremality of $\phi$
implies that $\omega$ is a multiplicity of some pure state on $\fB(\cK)$ i.e. $\omega(X)=\langle\eta,X\eta\rangle$ for some $\eta \in\cK$.
Let $x\in\cH$ be such that $Q=xx^*$. Then we have $\phi(X)=\langle \eta,X\eta\rangle xx^*=x\eta^*X\eta x^*=
(x\eta^*)X(x\eta^*)^*$, so (3) is fulfiled. Note that in the similar way one can show that 
$\phi(X)=(x\ov{\eta}^*)X^T(x\ov{\eta}^*)^*$.
\end{proof}
From the last result we immediately obtain the following characterization of non-decomposable extremal maps
\begin{corollary}\label{c:wn2}
If $\phi$ is such an extremal positive map that $\rank\, \phi(P)\geq 2$ for
some one-dimensional projection $P\in\fB(\cK)$ then $\phi$ is non-decomposable.
\end{corollary}
\begin{remark}\label{r:Choi}
Having Corollary \ref{c:wn2} one can easily prove that the Choi map $\phi$ given by (\ref{Choi}) is non-decomposable.
Indeed it is enough to calculate that in general for a one-dimensional projection $P$ the
map $\phi$ takes a value being an invertible matrix. It can be calculated that the only exceptions are
the following possible four values of $P$
$$\left[\begin{array}{ccc}
1/3&1/3&1/3\\1/3&1/3&1/3\\1/3&1/3&1/3
\end{array}\right],\quad
\left[\begin{array}{ccc}
1&0&0\\0&0&0\\0&0&0
\end{array}\right],\quad
\left[\begin{array}{ccc}
0&0&0\\0&1&0\\0&0&0
\end{array}\right],\quad
\left[\begin{array}{ccc}
0&0&0\\0&0&0\\0&0&1
\end{array}\right],
$$
for which the values of $\phi$ are matrices of rank equal to $2$.
\end{remark}

\section{General case}
We assume that $\phi:\fB(\cK)\to\fB(\cH)$ is a bounded positive map.
If $\phi$ is non-zero then we may find unit vectors $\xi\in\cK$, $x\in\cH$ and a positive number $\lambda$
such that
\begin{equation}
\phi(\xi\xi^*)x=\lambda x.\label{e:cond0}
\end{equation}
Let us fix such $\xi$, $x$ and $\lambda$.
Define two bounded operators $B,C:\cK\to\cH$ by
\begin{eqnarray}
B\eta&=&\lambda^{-1/2}\phi(\eta\xi^*)x,\label{e:B}\\
C\eta&=&\lambda^{-1/2}\phi(\xi\overline{\eta}^*)x\label{e:C}
\end{eqnarray}
where $\eta\in\cK$ and let $\psi$ and $\chi$ be maps from $\fB(\cK)$ into $\fB(\cH)$ determined by
\begin{eqnarray}
\psi(X)&=&BXB^*,\label{e:psi}\\
\chi(X)&=&CX^{\mathrm T}C^*\label{e:chi}
\end{eqnarray}
for $X\in\fB(\cK)$ (cf. \cite{Sem06}). Then we have the following
\begin{proposition}
Assume that $\phi(X)=AXA^*$ (resp. $\phi(X)=AX^{\mathrm T}A^*$) for $X\in\fB(\cK)$ where $A\in\fB(\cK,\cH)$ 
is some non-zero operator.
Let $\xi$, $x$ and $\lambda$ fulfil (\ref{e:cond0}). Take the operator $B$ as in (\ref{e:B}) 
(resp. $C$ as in (\ref{e:C})) and the map $\psi$ as in (\ref{e:psi}) (resp. $\chi$ as in (\ref{e:chi})).
Then $B=e^{it}A$ (resp. $C=e^{it}A$) for some $t\in\bR$ and $\psi=\phi$ (resp. $\chi=\phi$).
\end{proposition}
\begin{proof}
It follows from (\ref{e:cond0}) that $\langle A\xi,x\rangle A\xi=\lambda x$. Hence $A\xi=ax$ for some $a\in\bC$ 
such that $|a|=\lambda^{1/2}$.
Let $\eta\in\cK$. We calculate
$$B\eta=\lambda^{-1/2}\phi(\eta\xi^*)x=\lambda^{-1/2}A\eta\xi^*A^*x=\lambda^{-1/2}\langle A\xi,x\rangle A\eta=
e^{it}A\eta,$$
where $t\in\bR$ is such that $a=\lambda^{1/2}e^{-it}$. Consequently,
$\psi(X)=BXB^*=AXA^*=\phi(X)$.

If $\phi(X)=AX^{\mathrm T}A^*$ then we observe as above that $A\overline{\xi}=ax$ for some $a\in\bC$ of the form
$a=\lambda^{1/2}e^{-it}$ with $t\in\bR$, and for $\eta\in\cK$
\begin{eqnarray*}
C\eta&=&\lambda^{-1/2}\phi(\xi\overline{\eta}^*)x=\lambda^{-1/2}A(\xi\overline{\eta}^*)^{\mathrm T}A^*x=
\lambda^{-1/2}A\eta\overline{\xi}^*A^*x=\\
&=&\lambda^{-1/2}\langle A\overline{\xi},x\rangle A\eta=e^{it}A\eta.
\end{eqnarray*}
As previously, $\chi(X)=CX^{\mathrm T}C^*=AX^{\mathrm T}A^*=\phi(X)$ for $X\in\fB(\cK)$.
\end{proof}

Motivated by the above result we want to investigate if there are some relations between arbitrary map $\phi$ satisfying 
(\ref{e:cond0}) and maps $\psi$ and $\chi$. The answer for this question is contained in next result.
\begin{theorem}\label{p:psil}
Assume that $\phi:\fB(\cK)\to\fB(\cH)$ is an arbitrary given positive non-zero map.
Let $\xi$, $x$ and $\lambda$ fulfil (\ref{e:cond0}). Define the operator $B$ as in (\ref{e:B}) (resp. $C$ as in
(\ref{e:C})) and let $\psi(X)=BXB^*$ (resp. $\chi(X)=CX^{\mathrm T}C^*$) for $X\in\fB(\cK)$.
Then $\psi\leq\phi$ if and only if for any $\eta\in\cK$ and $y\in\cH$ the following
inequality holds
\begin{equation}\label{e:psil}
|\langle y,\phi(\eta\xi^*)x\rangle|^2\leq\langle x,\phi(\xi\xi^*)x\rangle\langle y,\phi(\eta\eta^*)y\rangle.
\end{equation}
Analogously, $\chi\leq\phi$ if and only if for any $\eta\in\cK$ and $y\in\cH$ the following
inequality holds
\begin{equation}\label{e:chil}
|\langle y,\phi(\xi\eta^*)x\rangle|^2\leq\langle x,\phi(\xi\xi^*)x\rangle\langle y,\phi(\eta\eta^*)y\rangle.
\end{equation}
\end{theorem}
\begin{proof}
The condition $\psi\leq\phi$ holds if and only if for any $\eta\in\cK$ and $y\in\cH$ 
\begin{equation}
\langle y,\psi(\eta\eta^*)y\rangle\leq\langle y,\phi(\eta\eta^*)y\rangle.
\end{equation}
The left hand side of the above inequality is equal to
\begin{equation}
\langle y,\psi(\eta\eta^*)y\rangle = \langle y,B\eta\eta^*B^*y\rangle =
|\langle y,B\eta\rangle|^2=\lambda^{-1}|\langle y,\phi(\eta\xi^*)x\rangle|^2.
\end{equation}
Taking into account that $\langle x,\phi(\xi\xi^*)x\rangle=\lambda$ (cf. (\ref{e:cond0})) we obtain (\ref{e:psil}). 
The second part of the proposition can be proved by similar arguments.
\end{proof}

Let us remind that if $\phi$ is a unital (i.e. such that $\phi(\jed)=\jed$) map then $\phi$ is called a {\it Schwarz map}
if $\phi(X^*X)\geq \phi(X^*)\phi(X)$ for any $X\in\fB(\cK)$. In the context of not necessarily unital maps we can
adopt the concept of {\it local complete positivity} introduced in \cite{St63}.
It is important to note that due to \cite[Theorem 7.4]{St63} a map $\phi$ is locally completely positive
if and only if there is a constant $\gamma>0$ such that
the inequality
\begin{equation}\label{e:CS}
(\gamma\phi)(X^*X)\geq (\gamma\phi)(X^*)(\gamma\phi)(X)
\end{equation}
is satisfied for all $X\in\fB(\cK)$.
Analogously, we will say that a map $\phi$ is {\it locally completely copositive} if for every $X\in\fB(\cK)$
\begin{equation}\label{e:coCS}
(\gamma\phi)(XX^*)\geq (\gamma\phi)(X^*)(\gamma\phi)(X)).
\end{equation}


Robertson asked in \cite{Rob83} if there exist a Schwarz map between $C^*$-algebras which is extreme as a 
positive unital map, but which is not $2$-positive. In the framework of non-unital maps we can ask the following
\begin{quest}
Are there any
extremal locally completely positive maps which are not $2$-positive?
\end{quest}
Robertson showed that if we consider maps acting from $M_2(\bC)$ into $M_2(\bC)$ the answer for his question is negative. 
We will extend this result for more cases in the sequel. 
Now, we show a little bit weaker but general result. 
\begin{theorem}\label{t:2pos}
Assume that $\phi$ is non-zero and extremal in the cone of all positive maps from $\fB(\cK)$ into $\fB(\cH)$. If 
$\phi$ is $2$-positive (resp. $2$-copositive) then it is completely positive (resp. completely copositive) map.
\end{theorem}
\begin{proof}
If $\phi$ is non-zero then there are $\xi$, $x$ and $\lambda$ which fulfil (\ref{e:cond0}). Let $\eta\in\cK$.
Consider the matrix $$\left[\begin{array}{cc}\xi\xi^*&\xi\eta^*\\\eta\xi^*&\eta\eta^*\end{array}\right].$$
One can easily
show that it is a positive element of $M_2(\fB(\cK))$. Hence, from $2$-positivity of $\phi$ we conclude that
the matrix $$\left[\begin{array}{cc}\phi(\xi\xi^*)&\phi(\xi\eta^*)\\\phi(\eta\xi^*)&\phi(\eta\eta^*)\end{array}\right]$$
is a positive element of $M_2(\fB(\cH))$. It implies that for any $y\in\cH$ we have
$$\left|\begin{array}{cc}\langle x,\phi(\xi\xi^*)x\rangle &\langle x,\phi(\xi\eta^*)y\rangle\\\langle y,\phi(\eta\xi^*)x\rangle
&\langle y,\phi(\eta\eta^*)y\rangle\end{array}\right|\geq 0.$$
So, from Proposition \ref{p:psil} (cf. (\ref{e:psil}))
we infer that $\psi\leq\phi$ where $\psi$ is the map defined in (\ref{e:psi}). But $\phi$ is extremal, so it is a positive multiplicity of $\psi$, hence it is completely positive map.
The "copositive" part of the theorem can be proved by similar arguments.
\end{proof}

\section{Functions with parallelogram identity}
Assume that $\phi:\fB(\cK)\to\fB(\cH)$ is a positive map. Then it is easy to check that for all $\xi,\eta\in\cK$
$$\phi((\xi+\eta)(\xi+\eta)^*)+\phi((\xi-\eta)(\xi-\eta)^*)=2\phi(\xi\xi^*)+2\phi(\eta\eta^*).$$
This is a motivation for the following
\begin{definition} 
Let $\cK$ and $\cH$ be (not necessarily finite dimensional) Hilbert spaces.
We say that a function $R:\cK\to\fB(\cH)$ fulfils the {\it parallelogram identity} if for all $\xi,\eta\in\cK$
\begin{equation}\label{e:paral}
R(\xi+\eta)+R(\xi-\eta)=2R(\xi)+2R(\eta).
\end{equation}

If additionally $R(\eta)\in\fB(\cH)^+$ for any $\eta\in\cK$ then we say that $R$ is a {\it positive function with parallelogram identity}.
\end{definition}

The purpose of this section is to show that under some additional assumptions on a function with parallelogram
identity $R$ it is possible to reconstruct a positive map $\phi$ such that $R(\eta)=\phi(\eta\eta^*)$ for
any $\eta\in\cK$.

Firstly, we characterize scalar positive functions with parallelogram identity.
\begin{lemma}
Let $\cK$ be a Hilbert space and assume that $\mu:\cK\to\bR$ is a continuous positive function with parallelogram identity such that for any $\xi\cK$
\begin{equation}\label{e:symm}
\mu(-\xi)=\mu(\xi),\quad\mu(i\xi)=\mu(\xi)
\end{equation}
Then there is a positive operator $M$ on $\cK$ such that for any $\xi\in\cK$
\begin{equation}
\mu(\xi)=\langle \xi,M\xi\rangle.
\end{equation}
\end{lemma}
\begin{proof}
We apply main arguments from \cite{JvN35}. For the readers convenience we give the full proof.
Firstly, define for $\xi,\eta\in\cK$
\begin{equation}\label{e:real}
(\xi,\eta)_{\bR}=\frac{1}{4}\left(\mu(\xi+\eta)-\mu(\eta-\xi)\right).
\end{equation}
It follows from (\ref{e:symm}) that 
\begin{equation}\label{e:symm1}
(\xi,\eta)_{\bR}=(\eta,\xi)_{\bR}
\end{equation}
for any
$\xi,\eta\in\cK$, so the form $(\cdot,\cdot)_{\bR}$ is symmetric.
Let $x,y,z$ be arbitrary elements of $\cK$.
If we substitute $\xi=y+z$ and $\eta=x$ then from (\ref{e:paral}) we obtain
\begin{equation}\label{e:par1}
\mu(y+z+x)+\mu(y+z-x)=2\mu(y+z)+2\mu(x).
\end{equation}
On the other hand, by taking $\xi=y-z$ and $\eta=x$ we get
\begin{equation}\label{e:par2}
\mu(y-z+x)+\mu(y-z-x)=2\mu(y-z)+2\mu(x).
\end{equation}
By subtracting (\ref{e:par2}) from (\ref{e:par1}) we obtain
$$
\mu(y+x+z)-\mu(y+x-z)+\mu(y-x+z)-\mu(y-x-z)=2[\mu(y+z)-\mu(y-z)]
$$
which in the context of (\ref{e:real}) is equivalent to
\begin{equation}\label{e:par3}
(z,y+x)_{\bR}+(z,y-x)_{\bR}=2(z,y)_{\bR}.
\end{equation}
By taking $y=z$ we have
\begin{equation}
(z,2y)_{\bR}=2(z,y)_{\bR}.
\end{equation}
Moreover, for any $y',x'\in\cK$ we can put in (\ref{e:par3}) $y=\frac{1}{2}(y'+x')$, $x=\frac{1}{2}(y'-x')$ to get
\begin{equation}\label{e:add}
(z,y')_{\bR}+(z,x')_{\bR} 
=(z,y'+x')_{\bR}.
\end{equation}
Taking (\ref{e:symm1}) into account we conclude that $(\cdot,\cdot)_{\bR}$ is additive in both coordinates.

Now, let $S=\{\alpha\in\bR:\,(\xi,\alpha\eta)_{\bR}=\alpha(\xi,\eta)_{\bR},\;\xi,\eta\in\cK\}$.
Obviously, $1\in S$. It follows from (\ref{e:add}) that $\alpha+\beta\in S$ for $\alpha,\beta\in S$.
Moreover, (\ref{e:real}) implies $(\xi,0)_\bR=0$ for any $\xi\in\cK$, so we infer that $-\alpha\in S$ for any
$\alpha\in S$. Hence $\bZ\subset S$. Now assume that $\alpha,\beta\in S$, $\beta\neq 0$. 
Then
$\beta(\xi,\alpha\beta^{-1}\eta)_{\bR}=(\xi,\alpha\eta)_{\bR}=\alpha(\xi,\eta)_{\bR}$.
Thus $\alpha\beta^{-1}\in S$, and $\bQ\subset S$. It follows from continuity of $\mu$ that $(\cdot,\cdot)_{\bR}$ is
continuous in both coordinates. Hence $S$ is a closed subset of $\bR$ and consequently $S=\bR$.

Now, define
\begin{equation}\label{e:form}
(\xi,\eta)=(\xi,\eta)_{\bR}-i(\xi,i\eta)_{\bR}.
\end{equation}
Additivity of $(\cdot,\cdot)_{\bR}$ implies additivity of $(\cdot,\cdot)$. Moreover, we observe that
$(\xi,\alpha\eta)=\alpha(\xi,\eta)$ for $\alpha\in\bR$. Now, for  $\alpha,\beta\in\bR$  we calculate
\begin{eqnarray*}
(\xi,(\alpha+i\beta)\eta)&=&(\xi,(\alpha+i\beta)\eta)_{\bR}-i(\xi,(i\alpha-\beta)\eta)_{\bR}\\
&=&\alpha(\xi,\eta)_{\bR}+\beta(\xi,i\eta)_{\bR}-i\alpha(\xi,i\eta)_{\bR}+i\beta(\xi,\eta)_{\bR}\\
&=&(\alpha+i\beta)(\xi,\eta)_{\bR}-i(\alpha+i\beta)(\xi,i\eta)_{\bR}\\
&=&(\alpha+i\beta)(\xi,\eta).
\end{eqnarray*}
Hence $(\cdot,\cdot)$ is linear with respect to the second variable.
Moreover, from the second equality in (\ref{e:symm}) we have $(i\xi,i\eta)_{\bR}=(\xi,\eta)_{\bR}$ for $\xi,\eta\in\cK$, so
\begin{eqnarray*}
(\xi,\eta)&=&(\xi,\eta)_{\bR}-i(\xi,i\eta)_{\bR}=(\eta,\xi)_{\bR}-i(i\eta,\xi)_{\bR}=\\
&=&(\eta,\xi)_{\bR}-i(\eta,-i\xi)_{\bR}=
(\eta,\xi)_{\bR}+i(\eta,i\xi)_{\bR}=\overline{(\eta,\xi)}.
\end{eqnarray*}
We proved that $(\cdot,\cdot)$ is a continuous skew-symmetric form on $\cK$. Thus there exists $M\in\fB(\cK)$ such that
$(\xi,\eta)=\langle \xi,M\eta\rangle$.

Now, observe that (\ref{e:paral}) for $\eta=0$ takes the form 
$\mu(\xi)+\mu(\xi)=2\mu(\xi)+2\mu(0)$
so we deduce that $\mu(0)=0$.
We check that
$
(\xi,\xi)_{\bR}=\frac{1}{4}\left[\mu(2\xi)-\mu(0)\right]=\mu(\xi),
$
and
$
(\xi,i\xi)_{\bR}=\frac{1}{4}\left[\mu((1+i)\xi)-\mu((i-1)\xi)\right]=
\frac{1}{4}\left[\mu((1+i)\xi)-\mu(i(1+i)\xi)\right]=0,
$
so we get $\mu(\xi)=(\xi,\xi)=\langle\xi,M\xi\rangle$. Finally, it follows from the positivity of $\mu$ that $M$ is positive.
\end{proof}

\begin{theorem}\label{t:par}
Let $\cK$ and $\cH$ be Hilbert spaces.
Assume that $R:\cK\to\fB(\cH)$ is a positive map with parallelogram identity such that
\begin{enumerate}
\item[(i)] for every $\eta\in\cK$.
\begin{equation}\label{e:Rjedn}
R(-\eta)=R(i\eta)=R(\eta),
\end{equation}
\item[(ii)]
for any $\xi,\eta\in\cK$ the map $\bR\ni\alpha\mapsto R(\eta+\alpha\eta)\in\fB(\cH)$ is continuous at zero,
\item[(iii)]
for any $\varepsilon>0$ there are $n\in\bN$, $\zeta_1,\ldots,\zeta_n\in\cK$ and $\delta>0$ such that for any 
$m\in\bN$, $\xi_1,\ldots,\xi_m,\eta_1,\ldots,\eta_m\in\cK$ the condition 
$\Vert\sum_{i=1}^m\langle\xi_i,\zeta_j\rangle\eta_i\Vert<\delta$ for $j=1,\ldots,n$ implies
$\Vert\sum_{i=1}^m(R(\eta_i+\xi_i)-R(\eta_i-\xi_i))\Vert<\varepsilon$.
\end{enumerate}
Then there is a positive map $\phi:\fB(\cK)\to\fB(\cH)$ such that 
$
\phi(\eta\eta^*)=R(\eta)
$
for each $\eta\in\cK$. Moreover, the map is uniquely determined and continuous with respect
to strong topology in the domain and uniform topology in the predomain.
\end{theorem}
\begin{proof}
As in the proof of the previous theorem we define a 
map $[\cdot,\cdot]:\cK\times\cK\to\fB(\cH)$ by the formula
\begin{equation}
[\xi,\eta]=[\xi,\eta]_{\bR}-i[\xi,i\eta]_{\bR},
\end{equation}
where
\begin{equation}
[\xi,\eta]_{\bR}=\frac{1}{4}\left[R(\eta+\xi)-R(\eta-\xi)\right].
\end{equation}
By the same arguments as in the previous proof we show that $[\cdot,\cdot]_\bR$ is symmetric and additive in both coordinates and
$[\xi,\alpha\eta]_\bR=\alpha[\xi,\eta]_\bR$ for any $\xi,\eta\in\cK$ and $\alpha\in\bQ$. It follows from the continuity assumption
that for any $\xi,\eta\in\cK$ the function $\bR\ni\alpha\mapsto[\xi,\alpha\eta]_\bR\in\fB(\cH)$ is continuous. As in the previous proof we conclude that
$[\xi,\alpha\eta]_\bR=\alpha[\xi,\eta]_\bR$ for any $\alpha\in\bR$. Thus we can also show that the map $[\cdot,\cdot]$ is skew-symmetric form with
the property that $[\xi,\xi]\geq 0$ for all $\xi\in\cK$.

Now, for any $\xi,\eta\in\cK$ define $\phi(\eta\xi^*)=[\xi,\eta]$. 
This definition can be extended by linearity onto the subspace $\fB_{\mathrm f}(\cK)$ of all finite dimensional operators on $\cK$ provided
we will show that $\phi$ is properly defined. 
In order to show that $\phi$ is properly
defined one should prove that for any $n\in\bN$, $\xi_1,\ldots,\xi_n,\eta_1,\ldots,\eta_n\in\cK$ the equality
$\sum_{i=1}^n\eta_i\xi_i^*=0$ implies $\sum_{i=1}^n[\xi_i,\eta_i]=0$. Assume firstly that $\eta_1,\ldots,\eta_n$
are linearly independent. Then for any $\zeta\in\cK$ we have $\sum_{i=1}^n\langle\xi_i,\zeta\rangle\eta_i=0$, and 
consequently $\langle\xi_i,\zeta\rangle=0$ for any $i=1,\ldots,n$ and $\zeta\in\cK$. This leads to the conclusion
that $\xi_i=0$ for each $i=1,\ldots,n$. If $\eta_1,\ldots,\eta_n$ are dependent then let us choose
a maximal linearly independent subsystem, say $\eta_,\ldots,\eta_k$, of the system $\eta_1,\ldots,\eta_n$.
Then for any $j=k+1,\ldots,n$ we have $\eta_j=\sum_{i=1}^k\alpha_{ij}\eta_i$ for some coefficients
$\alpha_{ij}$. Thus
$\sum_{i=1}^n\eta_i\xi_i^*=\sum_{i=1}^k\eta_i\left(\xi_i+\sum_{j=k+1}^n\ov{\alpha_{ij}}\xi_j\right)^*$ and
we get $\xi_i+\sum_{j=k+1}^n\ov{\alpha_{ij}}\xi_j=0$ for each $i=1,\ldots,k$. Now, we can calculate
$\sum_{i=1}^n[\xi_i,\eta_i]=\sum_{i=1}^k\left[\xi_i+\sum_{j=k+1}^n\ov{\alpha_{ij}}\xi_j,\eta_i\right]=0$.

Now, condition (iii) implies that $\phi$ is continuous on the subspace $\fB_{\mathrm f}(\cK)$
on $\cK$ with respect to the strong topology. But $\fB_{\mathrm f}(\cK)$ is strongly dense in $\fB(\cK)$ (cf. \cite{BR}),
so $\phi$ can be uniquely extended to the whole $\fB(\cK)$.
\end{proof}

\section{Structural results}
As it was mentioned if $\phi$ is a non-zero map then always there exists such a triple $\xi$, $x$ and $\lambda$ that
the condition (\ref{e:cond0}) holds. From now on we will assume much stronger condition.
For any unit vectors $\xi\in\cK$ and $x\in\cH$ we define
\begin{equation}\label{e:G}
G_{\xi,x}=\{\phi\in\fP(\fB(\cK),\cH):\,\mbox{$\phi(\xi\xi^*)=\lambda xx^*$ for some $\lambda\geq 0$}\}.
\end{equation}
One can easily observe that for every $\xi$ and $x$ the set $G_{\xi,x}$ is a face of the cone of all positive maps.
Our goal is to describe all extremal positive maps which lay in faces of the above form.
\begin{remark}
Note that there are extremal positive maps which are outside of any face $G_{\xi,x}$. 
The map defined in (\ref{Choi}) can serve as an example (cf. Remark \ref{r:Choi}).
\end{remark}
Firstly, we formulate some properties of positive maps which belong to some face $G_{\xi,x}$ for some unit vectors $\xi\in\cK$ and $x\in\cH$.
To this end we need the following
\begin{lemma}\label{l:form}
Let $\cH$ be a Hilbert space and fix some unit vector $x\in \cH$. 
Then for any $Y\in\fB(\cH)$ there are uniquely defined $\alpha\in\bC$, $u,v\in\cH$ and $Z\in\fB(\cH)$ such that
\begin{enumerate}
\item $\langle x,u\rangle=0=\langle x,v\rangle$,
\item $\langle x,Zy\rangle=0=\langle x,Z^*y\rangle$ for any $y\in\cH$,
\item $Y=\alpha xx^*+ux^*+xv^*+Z$.
\end{enumerate}

Moreover, $Y\geq 0$ if and only if $\alpha\geq 0$, $Z\geq 0$, $u=v$, and 
\begin{equation}\label{e:wyzn}
uu^*\leq\alpha Z.
\end{equation}
\end{lemma}
\begin{proof}
We define
\begin{eqnarray*}
\alpha & = & \langle x,Yx\rangle ,\\
u & = & (\jed-xx^*)Yx,\\
v & = & (\jed-xx^*)Y^*x,\\
Z & = & (\jed-xx^*)Y(\jed-xx^*).
\end{eqnarray*}
One can verify the properties (1)--(3) as well as the uniqueness of $\alpha,u,v,Z$.

Now, if $Y\geq 0$ then $Y$ is selfadjoint and
$$\alpha xx^*+ux^*+xv^*+Z=Y=Y^*=\overline{\alpha}xx^*+vx^*+xu^*+Z^*.$$ From the uniqueness it follows that $\alpha\in\bR$, $Z$ is selfadjoint and $u=v$.
Let $y\in \cH\ominus\bC x$. Since $Y\geq 0$ we have $\langle ax+y,Y(ax+y)\rangle\geq 0$
for any $a\in\bC$. 
\begin{equation}
\label{e:a}
\alpha |a|^2+2\Re\; a\langle y,u\rangle +\langle y,Zy\rangle\geq 0.
\end{equation} 
Considering the case $a=1$ and $y=0$ we show that $\alpha\geq 0$ while taking $a=0$ and $y\in\cH\ominus\bC x$
shows that $Z\geq 0$.
Let $t\in\bR$ be such that $e^{it}\langle y,u\rangle=|\langle y,u\rangle|$.
Puting $e^{it}a$ or $-e^{it}a$ instead of $a$ in (\ref{e:a}) we obtain
$$
\alpha a^2 +2|\langle y,u\rangle|a+\langle y,Zy\rangle\geq 0
$$
for every $a\in\bR$. This is equivalent to $|\langle y,u\rangle|^2\leq\alpha\langle y,Zy\rangle$
and this leads to (\ref{e:wyzn}).
\end{proof}

\begin{proposition}\label{p:xieta*}
Assume that $\phi\in G_{\xi,x}$, so 
\begin{equation}\label{e:cond1}
\phi(\xi\xi^*)=\lambda xx^*
\end{equation}  for some unit vectors $\xi\in\cK$, $x\in\cH$ and nonnegative constant $\lambda$. 
Then for any $\eta\in\cK$ there are $\beta\in\bC$ and $u,v\in\cH$ such that
$\langle u,x\rangle=0=\langle v,x\rangle$ and
\begin{equation}\label{e:etaxi*}
\phi(\eta\xi^*)=\beta xx^*+ux^*+xv^*.
\end{equation}
\end{proposition}
\begin{proof}
For any $y\in\cH$ consider a positive linear function $\omega_y(Y)=\langle y,Yy\rangle$, $Y\in\fB(\cH)$, on
the algebra $\fB(\cH)$. Since $\phi$ is positive, $\omega_y\circ\phi$ is a positive functional on
the algebra $\fB(\cK)$. But every positive functional is automatically completely positive. It implies
that for any $\eta\in\cK$ the complex matrix
$$\left[\begin{array}{cc}\langle y,\phi(\xi\xi^*)y\rangle & \langle y,\phi(\xi\eta^*)y\rangle \\
\langle y,\phi(\eta\xi^*)y\rangle & \langle y,\phi(\eta\eta^*)y\rangle\end{array}\right]$$
is positive, so applying (\ref{e:cond1}) we get
$|\langle y,\phi(\eta\xi^*)y\rangle|^2\leq\lambda|\langle y,x\rangle|^2\langle y,\phi(\eta\eta^*)y\rangle$.
This implies that $\langle y,\phi(\eta\xi^*)y\rangle=0$ 
for any $y$ which is orthogonal to $x$. By polarization formula we conclude that
$\langle y,\phi(\eta\xi^*)z\rangle=0$ for any $y,z\in\cH\ominus\bC x$, and finally
$(\jed-xx^*)\phi(\eta\xi^*)(\jed-xx^*)=0$.
We finish the proof by applying Lemma \ref{l:form}.
\end{proof}
Let $\cH_x=\cH\ominus\bC x$. Lemma \ref{l:form} and Proposition \ref{p:xieta*} imply that for any $\phi\in G_{\xi,x}$ 
there are
functions $\beta,\mu:\cK\to\bC$, $u,v,r:\cK\to\cH_x$ and $R:\cK\to\fB(\cH_x)$ such that
\begin{eqnarray}
\phi(\eta\xi^*)&=&\beta(\eta)xx^*+xv(\eta)^*+u(\eta)x^*,\label{e:phietaxi}\\
\phi(\xi\eta^*)&=&\overline{\beta(\eta)}xx^*+xu(\eta)^*+v(\eta)x^*,\label{e:phixieta}\\
\phi(\eta\eta^*)&=&\mu(\eta)+xr(\eta)^*+r(\eta)x^*+R(\eta).\label{e:phietaeta}
\end{eqnarray}
Let $\eta\in\cK$ be fixed. 
We will not write the arguments of the above functions when it will not cause a confusion. From 
positivity of $\phi$ it follows that $\phi(\eta\eta^*)\geq 0$, thus 
by Lemma \ref{l:form}
\begin{equation}
\mu\geq 0,\qquad
R\geq 0,\qquad
rr^*\leq R.
\end{equation}
We prove some properties of functions which appear in formulas (\ref{e:phietaxi}), (\ref{e:phixieta}) and (\ref{e:phietaeta}).
\begin{proposition}\label{p:char1}
Let $\eta\in\cK$. Then for any number $\sigma\geq 0$, vector $s\in\cH_x$ and operator $S\in\fB(\cH_x)$ such that
$\sigma>0$, $S\geq 0$, $\Tr S<\infty$ and $ss^*\leq\sigma S$ we have
\begin{equation}\label{e:klucz}
|\sigma\beta+\langle s,u\rangle +\langle v,s\rangle|^2\leq \sigma\lambda\left(\sigma\mu+2\Re\langle s,r\rangle +\Tr(SR)\right)
\end{equation}
\end{proposition}
\begin{proof}
Since $\phi$ is a positive map $\omega\circ\phi$ is a positive functional on $\fB(\cK)$ for every
positive normal functional $\omega$ on $\fB(\cH)$. But a positive functional is automatically a completely
positive map, thus 
$$\left[\begin{array}{cc}\omega\circ\phi(\xi\xi^*)&\omega\circ\phi(\xi\eta^*)\\
\omega\circ\phi(\eta\xi^*)&\omega\circ\phi(\eta\eta^*)\end{array}\right]$$
is a positive matrix and consequently
\begin{equation}\label{e:det}
|\omega\circ\phi(\eta\xi^*)|^2\leq\omega\circ\phi(\xi\xi^*)\,\cdot\,\omega\circ\phi(\eta\eta^*).
\end{equation}
Now, let $\sigma$, $s$ and $S$ fulfil the assumptions of the proposition and
$$\rho=\sigma xx^*+xs^*+sx^*+S.$$
It follows from Lemma \ref{l:form} that $\rho$ is a positive trace class operator on $\cH$. It determines positive normal functional
$\omega_\rho(Y)=\Tr(\rho Y)$, $Y\in\fB(\cH)$.
Let us calculate
\begin{eqnarray*}
\omega_\rho\circ\phi(\xi\xi^*)&=&\sigma\lambda,\\
\omega_\rho\circ\phi(\eta\eta^*)&=&\sigma\mu+2\Re\langle s,r\rangle +\Tr(SR),\\
\omega_\rho\circ\phi(\xi\eta^*)&=&\sigma\beta+\langle s,u\rangle+\langle v,s\rangle.
\end{eqnarray*}
If we substitute the above expressions into (\ref{e:det}) then we obtain (\ref{e:klucz}).
\end{proof}

\begin{proposition}
Let $\eta\in\cK$. Then 
\begin{equation}\label{e:beta}
|\beta|^2\leq\lambda\mu.
\end{equation}
and for any $y\in\cH_x$ the following inequalities hold
\begin{equation}\label{e:R}
|\langle y,u\rangle +\langle v,y\rangle|^2\leq \lambda\langle y,Ry\rangle
\end{equation}
\begin{equation}\label{e:re}
\left[\Re\langle y,\lambda r-\overline{\beta}u-\beta v\rangle\right]^2\leq (\lambda\mu-|\beta|^2)\left(\lambda\langle y,Ry\rangle - |\langle y,u\rangle +\langle v,y\rangle|^2\right).
\end{equation}
\end{proposition}
\begin{proof}
Let $\sigma$ be a positive number, $s=y$, and $S=\sigma^{-1}yy^*$. Then $\sigma,s,S$ fulfil the assumption of Proposition \ref{p:char1}, and inequality (\ref{e:klucz}) takes the form
$$
|\sigma\beta+\langle y,u\rangle +\langle v,y\rangle|^2\leq \lambda\mu\sigma^2+2\lambda\sigma\Re\langle y,r\rangle +\lambda\langle y,Ry\rangle .
$$
It can be rewritten in the form
$$ 
\sigma^2|\beta|^2+2\sigma\Re\left(\overline{\beta}(\langle y,u\rangle +\langle v,y\rangle)\right)+|\langle y,u\rangle +\langle v,y\rangle|^2\leq 
\lambda\mu\sigma^2+2\lambda\sigma\Re\langle y,r\rangle +\lambda\langle y,Ry\rangle
$$ 
and finally
$$(\lambda\mu-|\beta|^2)\sigma^2 +2\sigma\Re\langle y,\lambda r-\overline{\beta}u-\beta v\rangle+\lambda\langle y,Ry\rangle-|\langle y,u\rangle +\langle v,y\rangle|^2\geq 0.$$
Since the inequality holds for any $\sigma>0$ inequalities (\ref{e:beta}) and (\ref{e:R}) are evident.
By considering $-y$ instead of $y$ if necessary we conclude that the above inequality holds for every $\sigma\in\bR$.
But this statement is equivalent to inequality (\ref{e:re}).
\end{proof}

\begin{proposition}\label{p:psilw}
Let $\eta\in\cK$. The inequality (\ref{e:psil}) holds for any $y\in\cH$ if and only if 
\begin{equation}\label{e:wu}
(\lambda r-\overline{\beta}u)(\lambda r-\overline{\beta}u)^*\leq(\lambda\mu-|\beta|^2)(\lambda R-uu^*).
\end{equation}

Analogously, the inequality (\ref{e:chil}) is equivalent to 
\begin{equation}
(\lambda r-\beta v)(\lambda r-\beta v)^*\leq(\lambda\mu-|\beta|^2)( \lambda R-vv^*).
\end{equation}
\end{proposition}
\begin{proof} From (\ref{e:phietaxi}) we have
$$|\langle y,\phi(\eta\xi^*)x\rangle|^2 =|\langle y,\beta x+u\rangle|^2=\langle y,(\beta x+u)
(\beta x+u)^*y\rangle.$$
So, it follows from (\ref{e:phietaeta}) that (\ref{e:psil}) is equivalent to
$$(\beta x+u)
(\beta x+u)^*\leq \lambda(\mu xx^*+xr^*+rx^*+R).$$
The above inequality can be rewritten as
$$(\lambda\mu-|\beta|^2)xx^* +x(\lambda r-\overline{\beta} u)^*-(\lambda r-\overline{\beta} u)x^*+\lambda R-uu^*\geq 0.$$
By Lemma \ref{l:form} this is equivalent to (\ref{e:wu}).
The second part of the proposition can be proved by similar arguments.
\end{proof}

For any unit vectors $\eta\in\cK$ and $y\in\cH$ let us define
\begin{equation}\label{e:max_face}
F_{\eta,y}=\{\phi\in\fP(\fB(\cK),\cH):\,\phi(\eta\eta^*)y=0\}
\end{equation}
One can easily check that it is a face of the cone of all positive maps.
Let us recall that Kye (\cite{Kye96}) showed that each maximal face in the cone $\fP(\fB(\cK),\cH)$
is of the above form for some $\eta$ and $y$ provided that $\cK$ and $\cH$ are finite dimensional.
Observe that for any $\xi$ and $x$ we have
$$G_{\xi,x}=\bigcap_{y\bot x}F_{\xi,y}.$$
Moreover, $\phi\in G_{\xi,x}\cap F_{\xi,x}$ if and only if $\phi(\xi\xi^*)=0$.
\begin{theorem}\label{p:LCP}
Let $\cK$ and $\cH$ be arbitrary Hilbert spaces, and $\xi\in\cK$, $x\in\cH$ be unit vectors.
Assume that $\phi\in G_{\xi,x}\setminus F_{\xi,x}$ is locally completely positive. 
Then $\psi\leq\phi$ where $\psi$ is the completely
positive map defined in (\ref{e:psi}).
\end{theorem}
\begin{proof}
By the assumption $\phi(\xi\xi^*)=\lambda xx^*$ for some $\lambda>0$.
For $\eta\in\cK$ let $X=\eta\xi^*$. Then the inequality (\ref{e:CS}) (see also (\ref{e:phietaxi}) and (\ref{e:phixieta})) leads to 
$$\gamma^2(\overline{\beta} xx^*+xu^*+vx^*)(\beta xx^*+xv^*+ux^*)\leq \gamma \lambda xx^*$$
The left hand side of the above inequality is equal to
$$\gamma\left((|\beta|^2+\Vert u\Vert^2)xx^* + \overline{\beta} xv^* +  \beta vx^* + vv^*\right).$$
It is majorized by a multiplicity of the $1$-dimensional projection $xx^*$, so we conclude that $v=0$.

Let $y\in\cH_x$. Then inequality (\ref{e:re}) takes the form
\begin{equation}\label{e:reu}
\left[\Re\langle y,\lambda r-\overline{\beta}u\rangle\right]^2\leq
(\lambda\mu-|\beta|^2)(\lambda\langle y,Ry\rangle -|\langle y,u\rangle|^2).
\end{equation}
Let $t\in\bR$ be such a number that $e^{-it}\langle y,\lambda r-\overline{\beta}u\rangle=|\langle y,\lambda r-\overline{\beta}u\rangle|$. If we put $e^{it}y$ instead of $y$ in (\ref{e:reu}) then we obtain
$$
\langle y,(\lambda r-\overline{\beta}u)(\lambda r-\overline{\beta}u)^*y\rangle\leq
(\lambda\mu-|\beta|^2)\langle y,(\lambda R-uu^*)y\rangle.
$$
Since the above inequality is valid for any $y\in\cH_x$ the condition (\ref{e:wu}) is fulfiled.
Now, we take into account Propositions \ref{p:psilw} and \ref{p:psil} to conclude that $\psi\leq\phi$.
\end{proof}

Now we are ready to formulate results which give a partial answer for Robertson's question.
First of them establish negative answer in general case if we restrict ourselves to maps 
contained in $G_{\xi,x}\setminus F_{\xi,x}$.
\begin{theorem}\label{t:G-F}
Assume that a positive map $\phi:\fB(\cK)\to\fB(\cH)$ fulfils the following conditions
\begin{enumerate}
\item $\phi\in  G_{\xi,x}\setminus F_{\xi,x}$ for some unit vectors $\xi\in\cK$, $x\in\cH$,
\item $\phi$ is extremal in the cone of positive maps,
\item $\phi$ is locally completely positive.
\end{enumerate}
Then $\phi$ is of the form $\phi(X)=BXB^*$ for some bounded linear operator $B\in\fB(\cK,\cH)$.
\end{theorem}
\begin{proof}
It is immediate consequence of the previous theorem.
\end{proof}
Our next result establishes negative answer for Robertson's question in some special cases.
\begin{corollary}\label{t:G-F2}
Assume that $\cK$ is any finite dimensional Hilbert space and $\dim\cH=2$.
Then any locally completely positive map which is extremal in the cone of positive maps between
$\fB(\cK)$ and $\fB(\cH)$ is completely positive.
\end{corollary}
\begin{proof}
We may assume that $\phi$ is non-zero. Let $\cK_0=\{\eta\in\cK:\,\phi(\eta\eta^*)=0\}$.
It is a subspace of $\cK$. Indeed, it follows from Theorem \ref{t:par} that for any $\eta_1,\eta_2\in\cK_0$
and $\alpha_1,\alpha_2\in\bC$ we have 
\begin{eqnarray*}
0&=&2|\alpha_1|^2\phi(\eta_1\eta_1^*)+
2|\alpha_2|^2\phi(\eta_2\eta_2^*)\\
&=&\phi((\alpha_1\eta_1+\alpha_2\eta_2)(\alpha_1\eta_1+\alpha_2\eta_2)^*)+
\phi((\alpha_1\eta_1-\alpha_2\eta_2)(\alpha_1\eta_1-\alpha_2\eta_2)^*).
\end{eqnarray*}
Thus, in particular, $\phi((\alpha_1\eta_1+\alpha_2\eta_2)(\alpha_1\eta_1+\alpha_2\eta_2)^*)=0$.
Let $P$ be the projection onto $\cK_0$ and $Q=\jed -P$. Then for any $X\in\fB(\cK)$ we have
$\phi(X)=\phi(QXQ)$. Let $\phi':Q\fB(\cK)Q\to\fB(\cH)$ be the compression of $\phi$ onto the algebra
$Q\fB(\cK)Q$. We show that the map $\phi'$ is extremal in the cone of all positive maps between 
$Q\fB(\cK)Q$ and
$\fB(\cH)$. Assume $\rho'\leq\phi'$ for some positive map $\rho':Q\fB(\cK)Q\to\fB(\cH)$ and define
$\rho:\fB(\cK)\to\fB(\cH)$ by $\rho(X)=\rho'(QXQ)$ for $X\in\fB(\cK)$. Then $\rho=\alpha\phi$ for some $\alpha\geq 0$ because $\phi$ is extremal. But this implies $\rho'=\alpha\phi'$.

Since $\phi'$ is extremal it must contain in some maximal face of the cone of all positive maps between
$Q\fB(\cK)Q$ and $\fB(\cH)$.
By the result of Kye it follows that there are unit vectors $\xi\in\cK\ominus\cK_0$ and $y\in\cH$
such that $\phi'(\xi\xi^*)y=0$. The condition $\dim\cH=2$ implies that 
$\phi'(\xi\xi^*)=\lambda xx^*$ where $x\in\cH$ is a unit vector such that $x\bot y$. From the definition of $\phi'$ it follows that $\phi'(\xi\xi^*)$ is non-zero, so $\lambda>0$.
Thus we proved $\phi(\xi\xi^*)=\phi(Q\xi\xi^*Q)=\psi'(\xi\xi^*)=\lambda xx^*$, and consequently
$\phi\in G_{\xi,x}\setminus F_{\xi,x}$. The rest follows from Theorem \ref{t:G-F}.
\end{proof}

\begin{remark}
Let us note that Theorem \ref{t:G-F} gives also a partial solution of the problem described
in Remark \ref{rem1}. That problem can be reformulated as follows: {\it Is there an extremal map
in $G_{\xi,x}\setminus F_{\xi,x}$ such that it is not a positive rank 1 nonincreasing map?}
Theorem \ref{t:G-F} provides a negative answer to this question if we restrict our considerations to the class
of locally completely positive extremal maps.
\end{remark}

\vspace{3mm} \textit{Acknowledgments.} The author want to thank to Adam Majewski and Louis E. Labuschagne for 
inspiration and fruitful discussions.

\end{document}